\documentclass{amsart}
\usepackage{amssymb,amscd}

\newtheorem*{theorem*}{Theorem}
\newtheorem{theorem}{Theorem}
\newtheorem{lemma}{Lemma}

\newtheorem{corollary}{Corollary}
\newtheorem*{corollary*}{Corollary}
\theoremstyle{remark}
\newtheorem*{remark}{Remark}

\theoremstyle{definition}
\newtheorem{definition}{Definition}
\newtheorem{example}{Example}

\newcommand{\thref}[1]{Theorem~{\rm\ref{#1}}}
\newcommand{\leref}[1]{Lemma~{\rm\ref{#1}}}
\newcommand{\seref}[1]{Section~{\rm\ref{#1}}}

\newcommand{\st}{\; | \;}                     
\newcommand{\ttt}{\otimes}                    

\newcommand{\<}{\langle}
\renewcommand{\>}{\rangle}

\newcommand{\C}{\mathbb{C}}       
\newcommand{\Z}{\mathbb{Z}}       
\newcommand{\R}{\mathbb{R}}       
\newcommand{\HH}{\mathbb{H}}      
\newcommand{\T}{\mathbb{T}}

\newcommand{\al}{\alpha}

\newcommand{\de}{\delta}

\newcommand{\la}{\lambda}

\newcommand{\om}{\omega}

\newcommand{\g}{\mathfrak{g}}
\newcommand{\aaa}{\mathfrak{a}}

\newcommand{\liet}{\mathfrak{t}}

\newcommand{\sll}{\mathfrak{sl}}

\newcommand{\h}{\mathfrak{h}}
\newcommand{\SU}{\mathrm{SU}}
\newcommand{\U}{\mathrm{U}}
\newcommand{\SO}{\mathrm{SO}}
\renewcommand{\O}{\mathrm{O}}
\newcommand{\Sp}{\mathrm{Sp}}
\newcommand{\Spin}{\mathrm{Spin}}

\newcommand{\GL}{\mathrm{GL}}

\DeclareMathOperator{\Lie}{Lie}
\DeclareMathOperator{\Ad}{Ad}

\DeclareMathOperator{\tr}{tr}
\DeclareMathOperator{\id}{id}
\DeclareMathOperator{\Aut}{Aut}
\DeclareMathOperator{\End}{End}
\DeclareMathOperator{\Hom}{Hom}
\DeclareMathOperator{\Mat}{Mat}
\DeclareMathOperator{\rk}{rk}
\DeclareMathOperator{\diag}{diag}
\DeclareMathOperator{\mult}{mult}
\DeclareMathOperator{\Res}{Res}

\begin{document}

\title{Compact groups and their representations}

\author{Alexandre Kirillov}
\author{Alexander Kirillov Jr}

\address{A.K.: Department of Mathematics\\
  University of Pennsylvania\\
  Philadelphia, PA 19104-6395, USA\\
  e-mail: kirillov@math.upenn.edu}
  \address{A. K. Jr:  Department of Mathematics\\
  SUNY at Stony Brook\\
  Stony Brook, NY 11794--3651, USA\\
  e-mail:kirillov@math.sunysb.edu}

\date{May 28, 2005}

\begin{abstract}
This is an overview article on compact Lie groups and their 
representations, written for the {\em Encyclopedia of Mathematical 
Physics} to be published by Elsevier.
\end{abstract}

\maketitle


In this article we describe the structure and representation theory of
compact Lie groups.  Throughout the article, $G$ is a compact real Lie 
group with Lie algebra $\g$. Unless otherwise stated, $G$ is
assumed to be connected. The word ``group'' will always mean a 
``Lie group'' and the word ``subgroup'' will mean a closed Lie 
subgroup. The notation $\Lie(H)$ stands for the Lie algebra of a Lie 
group $H$. We assume that the reader is familiar with the basic facts
of the theory of Lie groups and Lie algebras, which can be found in
the article {\em Lie groups: general theory} in this Encyclopedia, or
in  the books listed in the bibliography. 

\section{Examples of compact Lie groups} \label{s:examples}

Examples of compact groups include
\begin{itemize}
    \item Finite groups
    \item  Quotient groups $\T^n=\R^n/\Z^n$, or more
        generally, $V/L$, where $V$ is a finite-dimensional real
        vector space and $L$ is a lattice in $V$, i.e. a discrete
        subgroup generated by some basis in $V$.
        Groups of this type are called {\bf tori}; it is known that
        every commutative connected compact group is a {\bf torus}.

    \item Unitary groups $\U(n)$ and special unitary
        groups $\SU(n)$, $n\geq 2$.

    \item  Orthogonal groups $\O(n)$ and $\SO(n)$, $n\ge 3$.
  
    \item  The groups $U(n,\HH)$,  $n\ge 1$, of unitary
      quaternionic  transformations, which are isomorphic to
      $\Sp(n)\,:=\Sp(n,\C)\cap \SU(2n)$. 
\end{itemize}

The groups $\O(n)$  have two connected components, one of which
is $\SO(n)$. 

The groups $\SU(n)$ and $\Sp(n)$ are connected and simply-connected. 

The groups $\SO(n)$ are connected but not simply-connected: for $n\ge
3$, the fundamental group of $\SO(n)$ is $\Z_2$. The universal cover 
of $\SO(n)$ is a simply-connected compact Lie group denoted by
$\Spin(n)$. For small $n$, we have isomorphisms:
$\Spin(3)\simeq \SU(2)$,\quad $\Spin(4)\simeq \SU(2)\times
\SU(2)$,\quad
$\Spin(5)\simeq \Sp(4),\quad \Spin(6)\simeq \SU(4)$.

\section{Relation to semisimple Lie algebras and Lie groups}
\label{s:semisimple}

\subsection{Reductive groups}

A Lie algebra $\g$ is called 

\begin{itemize}
\item  {\bf simple} if it is non-abelian and has no ideals
  different from $\{0\}$ and $\g$ itself. 

\item  {\bf semisimple} if it is a direct sum of simple ideals.

\item  {\bf reductive} if it is a direct sum of  semisimple and 
   commutative ideals.
\end{itemize}
We call a connected Lie group $G$ {\bf simple} or {\bf semisimple} if 
$\Lie(G)$ has this property.

\begin{theorem}\label{t:reductive}
 Let $G$ be a connected compact Lie group and 
  $\g=\Lie(G)$. Then 
  \begin{enumerate}
    \item The Lie algebra $\g=\Lie(G)$ is reductive: $\g=\aaa\oplus
    \g'$, where $\aaa$ is abelian and $\g'=[\g,\,\g]$ is semisimple.

  \item The group $G$ can be written in the form $G=(A \times
    K)/Z$, where $A$ is a torus, $K$ is a connected,
    simply-connected compact semisimple Lie group, and
    $Z$ is a finite central subgroup in $A\times K$.

  \item If $G$ is simply-connected, it is a product of simple compact
    Lie groups.
  \end{enumerate} 
 \end{theorem}

 The proof of this results  is based on the
fact that the Killing form of $\g$ is negative semidefinite.

\begin{example} The group $\U(n)$ contains as the center
  the subgroup $C$ of scalar matrices. The quotient group $\U(n)/C$ is
  simple and isomorphic to $\SU(n)/\Z_n$. The presentation of
  \thref{t:reductive} in this case is 
  $$
  \U(n)=\bigl(\T^1\times \SU(n)\bigr)/\Z_n 
  =\bigl(C\times \SU(n)\bigr)/\bigl(C\cap \SU(n)\bigr).
  $$
  For the group $\SO(4)$ the presentation is
  $\bigl(\SU(2)\times\SU(2)\bigr)/\{\pm (1\times 1)\}$.
\end{example}

This theorem effectively reduces the study of the structure of
connected compact groups to the study of simply-connected compact
simple Lie groups.

\subsection{Complexification of a compact Lie group}

Recall that for a real Lie algebra $\g$, its {\bf complexification}
is $\g_\C=\g\otimes \C$ with obvious commutator. It is also well-known
that $\g_\C$ is semisimple or reductive iff $\g$ is semisimple or
reductive respectively. There is a subtlety in the case of simple 
algebras: it is possible that a real Lie algebra is simple, but its 
complexification $\g_\C$ is only semisimple.
However, this problem never arises for 
Lie algebras of compact groups: if $\g$ is a Lie algebra of a
real compact Lie group, then $\g$ is simple if and only if $\g_\C$ is
simple.

The notion of complexification for Lie groups is more delicate.

 \begin{definition}\label{d:complexification}
   Let $G$ be a connected real Lie group with Lie algebra $\g$.
  A {\bf complexification} of $G$ is a connected {\bf complex Lie
  group} $G_\C$ \textup{(}i.e. a complex manifold with a structure
  of a Lie  group such that group multiplication is given by a complex
  analytic map $G_\C \times G_\C\to G_\C$\textup{)} which contains
  $G$ as a closed subgroup, and such that $\Lie(G_\C)=\g_\C$. 
  In this case, we will also say that $G$ is a {\bf real form} of
  $G_\C$.
 \end{definition}

 It is not obvious why such a complexification exists at all; in fact,
 for arbitrary real group it may not exist. However, for compact
 groups we do have the following theorem.

\begin{theorem} \label{t:complexification}
   Let $G$ be a connected compact  Lie group. Then it has a
   unique complexification $G_\C\supset G$. Moreover, the following
   properties hold:
  \begin{enumerate}
     \item The inclusion $G\subset G_\C$ is a homotopy equivalence. In
       particular, $\pi_1(G)=\pi_1(G_\C)$ and the quotient space
       $G_\C/G$ is contractible.

   \item Every finite-dimensional representation of $G$ can be
     uniquely extended to a complex analytic representation of $G_\C$.
  \end{enumerate}
\end{theorem}

Since a Lie algebra of a compact Lie group $G$ is reductive, we see
that $G_\C$ must be reductive; if $G$ is semisimple or simple, then
so is $G_\C$. The natural question is whether every complex reductive
group can be obtained in this way. The following theorem gives a
partial  answer.

 \begin{theorem}  \label{t:compact_form}
  Every connected complex semisimple Lie group $H$  has a compact
  real form: there is a compact real subgroup $K\subset H$ such that
  $H=K_\C$.  Moreover, such a compact real form is unique up to
  conjugation.
\end{theorem}

 \begin{example}
   The  unitary group  $\U(n)$ is a compact real form of the group
   $\GL(n,\C)$.

    The  orthogonal  group  $\SO(n)$ is a compact real form of the
     group  $\SO(n,\C)$.

   The group $\Sp(n)$ is a compact real form of the group $\Sp(n,\C)$.

   The universal cover of $\GL(n,\C)$ has no compact real form.
\end{example}

 These results have a number of important applications. For example, 
they show that study of representations of a semisimple complex
group $H$ can be replaced by the study of representations of its
compact form; in particular, every representation is completely
reducible (this argument is known as Weyl's unitary trick).

\subsection{Classification of simple compact Lie groups}
\thref{t:reductive} essentially reduces such classification 
to classification of simply-connected simple compact groups, and
\thref{t:complexification}, \thref{t:compact_form}  
reduce it to the classification of simple complex Lie algebras. 
Since the latter is well-known, we get the following result.

\begin{theorem} \label{t:classification}
Let $G$ be a connected, simply-connected simple compact Lie group. 
Then $\g_\C$ must be a simple complex Lie
algebra and thus can be  described by a Dynkin diagram of one the
following types: $A_n$, $B_n$, $C_n$, $D_n$, $E_6$,
$E_7$, $E_8$, $F_4$, $G_2$.

 Conversely, for each Dynkin diagram in the above list, there exists 
a unique up to isomorphism simply-connected simple compact Lie group 
whose Lie algebra is described by this Dynkin diagram.
\end{theorem}

 For types $A_n,\dots, D_n$, the corresponding compact Lie groups are
 well-known classical groups shown in the table below:

\medskip
 \begin{tabular}{|c|c|c|c|}
 \hline
 $A_n, \, n\ge 1$ & $B_n,\, n\ge 2$  & $C_n, n\ge 3$ & $D_n,\,
n\ge 4 $\\
 \hline
 $\SU(n+1)$ & $\Spin(2n+1)$  & $\Sp(n)$& $\Spin(2n)$\\
 \hline
 \end{tabular}

 \medskip
 The restrictions on $n$ in this table are made to avoid repetitions
which appear for small values of $n$. Namely, $A_1=B_1=C_1$, which
gives $\SU(2)=\Spin(3)=\Sp(1)$; \quad $D_2=A_1\cup A_1$, which gives 
$\Spin(4)=\SU(2)\times \SU(2)$; $B_2=C_2$, which gives $\SO(5)=\Sp(4)$
and $A_3=D_3$, which gives $\SU(4)=\Spin(6).$ 
Other than that, all entries are distinct.

 Exceptional groups $E_6,\dots, G_2$ also admit explicit geometric and
algebraic descriptions which are related to the exceptional
non-associative algebra $\mathbb{O}$ of so-called octonions (or Cayley
numbers). For example, the compact group of type $G_2$ can be
defined as a subgroup of $\SO(7)$ which preserves an almost complex
structure on $S^6$. It can also be described  as the subgroup of
$\GL(7,\,\R)$ which preserves one quadratic and one cubic form, or,
finally, as a group of all automorphisms of $\mathbb{O}$.

\section{Maximal tori} \label{s:tori}
\subsection{Main properties}
In this section, $G$ is a compact connected Lie group.
\begin{definition}\label{d:max_torus}
  A {\bf maximal torus} in $G$ is a maximal connected commutative
  subgroup $T\subset G$.
\end{definition}

 The following theorem lists the main properties of maximal tori.

\begin{theorem}\label{t:tori}
  \par\indent
  \begin{enumerate}
  \item For every element $g\in G$, there exists a 
    maximal torus $T\ni g$.
  \item Any two maximal tori in $G$ are conjugate.
  \item  If $g\in G$ commutes with all elements of a maximal
    torus $T$, then $g\in T$.
  \item A connected subgroup $H\subset G$ is a maximal torus iff
    the Lie algebra $\Lie(H)$ is a maximal abelian subalgebra in $\Lie(G)$.
\end{enumerate}
\end{theorem}

\begin{example} Let $G=\U(n)$. Then the set $T$ of diagonal unitary
  matrices is a maximal torus in $G$; moreover, every maximal torus
  is of this form after a suitable unitary change of basis. In
  particular, this implies that  every element in $G$ is conjugate to
  a diagonal matrix.
\end{example}

\begin{example} Let $G=\SO(3)$. Then the set $D$ of diagonal
matrices is a maximal commutative subgroup in $G$, but not a torus.
Here $D$ consists of 4 elements and is not connected. 
\end{example}

\subsection{Maximal tori and Cartan subalgebras}
The study of maximal tori in compact Lie groups is closely related
to the study of Cartan subalgebras in reductive complex Lie algebras 
(see \cite{serre} for a definition of a Cartan subalgebra). 

\begin{theorem}\label{t:tori-cartan}
  Let $G$ be a connected compact  Lie group with Lie
  algebra $\g$, and let $T\subset G$ be a maximal torus in $G$,
  $\liet=\Lie(T)\subset \g$. Let $\g_\C$, $G_\C$ be the
complexification  of $\g$, $G$ as in \thref{t:complexification}. 
  
  Let $\h=\liet_\C\subset \g_\C$. Then $\h$ is a Cartan subalgebra in
$\g_\C$.  Conversely, every  Cartan subalgebra in $G_\C$ can be
obtained as   $\liet_\C$ for some  maximal torus $T\subset G$.
  
  \end{theorem}
  
  This allows us to use results about Cartan subalgebras (such as root
decomposition, properties of root systems, etc.) when studying compact
Lie groups.  From now on, we will denote by $R\subset \h^*$ the  root
system of $\g_\C$; it can be shown that in fact, $R\subset i\liet^*$.
We will also use notation $\al^\vee\subset i\liet$ for the coroot
corresponding to a root $\al\in R$. Definitions of these and related
notions can be found in \cite{serre} or in the article {\em Lie
groups: general theory} in this encyclopedia. 
%
%
%

\subsection{Weights and roots}
Let $G$ be semisimple. Recall that  the root lattice $Q\subset
i\liet^*$ is the abelian group generated by roots $\al\in R$, and let
the coroot lattice  $Q^\vee\subset i\liet$
be the abelian group generated by coroots $\al^\vee$, $\al\in R$. 
Define also the weight and coweight lattices by 
\begin{align*}
P&=\{\la\st \<\al^\vee,\la\>\in \Z\quad\forall\al\in R\}\subset
i\liet^*\\
P^\vee&=\{t\st \<t,\al\>\in \Z\quad\forall\al\in R\}\subset i\liet,
\end{align*}
where $\langle\,\cdot ,\cdot  \,\rangle$ is the pairing between
$\liet$ and the dual vector space $\liet^*$. 

It follows from the definition of root system that we have inclusions 
\begin{equation}\label{e:PQ}
\begin{aligned}
Q \subset &P \subset i\liet^*\\
Q^\vee \subset &P^\vee \subset i\liet
\end{aligned}
\end{equation}

Both $P,Q$ are lattices in $i\liet^*$; thus, the index $(P:Q)$ is
finite. It can be computed explicitly: if $\al_i$ is a basis of the
root system,  then
the fundamental weights $\om_i$ defined by 
\begin{equation}\label{e:fund_weights}
\<\al_i^\vee,\om_j\>=\de_{ij}
\end{equation}
form a basis of $P$. The simple roots $\al_i$ are related to
fundamental weights $\om_j$ by the Cartan matrix $A$: $\al_i=\sum
A_{ij}\om_j$. Therefore, $(P:Q)=(P^\vee:Q^\vee)=|\det A|$. 

Definitions of  $P$, $Q$, $P^\vee$, $Q^\vee$
 if $\g$ also make sense when $\g$ reductive but not semisimple.
However, in this case they
are no longer lattices: $\rk Q<\dim\liet^*$, and $P$ is not discrete.

We can now give more precise information about the structure of the
maximal torus. 
\begin{lemma} \label{l:L} 
  Let $T$ be a compact connected commutative 
  Lie group, and  $\liet=\Lie(T)$  its Lie algebra. Then the
  exponential map is surjective and preimage of unit is a lattice
  $L\subset \liet$. There is an isomorphism of Lie groups
   $$
   \exp\colon  \liet/L\to T.
   $$
  In particular, $T\simeq \R^r/\Z^r=\T^r$, $r=\dim\, T$.
\end{lemma}

Let $X(T)\subset i\liet^*$ the lattice dual to $\frac 1
{2\pi i} L$:
\begin{equation}\label{e:X}
X(T)=\{\la\in i\liet^*\st \<\la,l\>\in 2\pi i \Z \quad \forall l\in
L\}.
\end{equation}
 It is  called the 
{\bf character lattice} for $T$ (see \seref{s:examples_reps}). 

\begin{theorem} Let $G$ be a compact connected  Lie group, 
and let $T\subset G$ be a maximal torus in $G$.

Then $Q\subset X(T)\subset P$. Moreover, the group $G$ is uniquely 
determined by the Lie algebra $\g$ and the lattice $X(T)\in \liet^*$ 
which can be any lattice between $Q$ and $P$.
\end{theorem}
\begin{corollary*}
For a given complex semisimple Lie algebra $\aaa$, there are only
finitely many \textup{(}up to isomorphism\textup{)} compact connected 
Lie groups $G$ with $\g_\C=\aaa$. 

The largest of them is the simply-connected group, for which $T=
\liet/2\pi i Q^\vee,\quad X(T)=P$; the smallest is the so-called 
{\bf adjoint group}, for which $T= \liet/2\pi i\,P^\vee,\quad
X(T)=Q$. 
\end{corollary*}

\begin{example}
Let $G=\U(n)$. Then $i\liet=\{\text{real diagonal  matrices}\}$.
Choosing the standard basis of matrix units in $i\liet$, we identify
$i\liet\simeq \R^n$, which also allows us to identify $i\liet^*\simeq
\R^n$. Under this identification, 
\begin{align*}
Q&=\{(\la_1,\dots, \la_n)\st \la_i\in\Z,\ \sum \la_i=0\},\\
P&=\{(\la_1,\dots, \la_n)\st \la_i\in\R,\ \la_i-\la_j\in \Z\},\\
X(T)&=\Z^n.
\end{align*}
Note that $Q,P$ are not lattices: $Q\simeq \Z^{n-1}$, $P\simeq
\R\times \Z^{n-1}$. 

Now let $G=\SU(n)$. Then $i\liet^*=\R^n/\R\cdot (1,\dots, 1)$, and
$Q,P$ are the images of $Q,P$ for $G=\U(n)$ in this quotient. In this
quotient they are lattices, and $(P:Q)=n$. The character lattice in
this case is $X(T)=P$, since $\SU(n)$ is simply-connected. The
adjoint group is $\mathrm{PSU}(n)=\SU(n)/C$, where $C=\{\la\cdot
\id\st\la^n=1\}$ is the center of $\SU(n)$. 

\end{example}

\subsection{Weyl group}

Let us fix a maximal torus $T\subset G$. Let $N(T)\subset G$ be the 
normalizer of $T$ in $G$: $N(T)=\{g\in G\bigm | gTg^{-1}=T\}$. For
any $g\in N(T)$ the transformation $A(g)\colon t\mapsto gtg^{-1}$ is
an automorphism of $T$. According to \thref{t:tori},
this automorphism is trivial iff $g\in T$. So in fact, it is the
quotient group $N(T)/T$ which acts on $T$.

\begin{definition}\label{d:weyl}
 The group $W=N(T)/T$ is called the {\bf Weyl group} of $G$. 
\end{definition}

Since the Weyl group acts faithfully on $\liet$ and $\liet^*$, it is 
common to consider $W$ as a subgroup in $\GL(\liet^*)$. 

It is known that $W$ is finite. 

The Weyl group can also be defined in terms of Lie algebra $\g$ and
its complexification $\g_{\C}$.  

 \begin{theorem}\label{t:weyl} 
 The Weyl group coincides with the subgroup in $\GL(i\liet^*)$
generated  
by reflections $s_\al\colon x\mapsto x-\frac {2
  (\al,x)}{(\al,\,\al)}$, $\al\in R$.
Here $(\, , \,)$ is the bilinear form on $\h^*$ induced by the Killing
form of $\g$. 
 \end{theorem}

 \begin{theorem}\label{t:T/W}
\par\indent
\begin{enumerate}
\item Two elements $t_1, t_2\in T$ are conjugate in $G$ 
iff $t_2=w(t_1)$ for some $w\in W$.
\item There exists a natural homeomorphism of quotient spaces $G/\Ad
 G\simeq T/W$, where $\Ad G$ stands for action of $G$ on itself by
 conjugation. \textup{(}Note, however, that these quotient spaces are
  not manifolds: they have singularities.\textup{)}

 \item Let us call a function $f$ on $G$ central if
$f(hgh^{-1})=f(g)$ for any $g,h\in G$. Then the restriction map
gives an isomorphism
$$\{\text{continuous central functions on $G$}\}
 \simeq \{\text{$W$-invariant continuous functions on $T$}\}$$
\end{enumerate}
\end{theorem}
 
\begin{example}
  Let $G=\U(n)$. The set of  diagonal   unitary matrices is a maximal
torus, and the Weyl group is the symmetric group $S_n$
  acting on diagonal matrices by permutations of entries. In this
  case, \thref{t:T/W} shows that if $f(U)$ is a central function of a
  unitary matrix, then $f(U)=\tilde f(\lambda_1,\dots,  \lambda_n)$,
  where $\lambda_i$ are eigenvalues of $U$ and $\tilde f$ is a
  symmetric function in $n$ variables.
\end{example}

\section{Representations of compact groups} \label{s:representations}
\subsection{Basic notions}

By a {\bf representation} of $G$ we understand a pair $(\pi,\,V)$
where $V$ is a complex vector space and  $\pi$ is a continuous
homomorphism $G\to \Aut(V)$. This notation is often shortened to $\pi$
or $V$. In this article we only consider finite-dimensional
representations;
in this case, the  homomorphism $\pi$ is automatically  smooth and
even real-analytic.
 
We associate to any f.d. representation $(\pi,\,V)$ of $G$ the
representation $(\pi_*,\,V)$ of the Lie algebra $\g=\Lie(G)$ which is
just the derivative of the map $\pi\colon  G\to \Aut V$ at the unit
point $e\in G$. In terms of the exponential map we have the following
commutative diagram
$$\begin{CD}
         G@>\pi>> \Aut\,V\\
       @A\exp AA @AA\exp A\\
         \g@ >\pi_*>>\End\,V.
\end{CD}
$$

Choosing a basis in $V$, we can write the operators $\pi(g)$ and 
$\pi_*(X)$ in matrix form and consider $\pi$ and $\pi_*$ as matrix
valued  functions on $G$ and $\g$. The diagram above means that
\begin{equation}\label{e:1}
  \pi(\exp\,X)= e^{\pi_*(X)}.
\end{equation}
Recall that if $G$ is connected, simply-connected, then every
representation of $\g$ can be uniquely lifted to a representation of
$G$. Thus, classification of representations of
connected simply-connected Lie groups is equivalent to the
classification of representations of Lie algebras.

Let $(\pi_1,\,V_1)$ and $(\pi_2,\,V_2)$ be two representation of the
same group $G$. An operator $A\in \Hom(V_1,\,V_2)$ is called an {\bf 
intertwining operator}, or simply an {\bf  intertwiner}, if
$A\circ\pi_1(g)=\pi_2(g)\circ A$ for all $g\in G$. Two
representations are called {\bf  equivalent} if
they admit an invertible intertwiner. In this case, using an
appropriate choice of bases, we can write $\pi_1$ and $\pi_2$ by the
same matrix-valued function.

Let $(\pi,\,V)$ be a representation of $G$. If all operators $\pi(g),
\,g\in G,$ preserve a subspace $V_1\subset V$, then the restrictions
$\pi_1(g)= \pi(g)|_{V_1}$ define a {\bf  subrepresentation}
$(\pi_1,\,V_1)$ of $(\pi,\,V)$. In this case,  the
quotient space $V_2=V/V_1$ also has a canonical structure of a
representation, called the {\bf quotient representation}.

A representation $(\pi,\,V)$ is called {\bf reducible} if it has
a non-trivial (different from $V$ and $\{0\}$) subrepresentation.
Otherwise it is called {\bf irreducible}.


We call representation $(\pi,\,V)$ {\bf  unitary} if $V$ is a Hilbert
space and all operators $\pi(g),\,g\in G,$ are unitary, i.e. given by
unitary matrices in any orthonormal basis. We use a short term
``unirrep'' for a ``unitary irreducible representation''.

\subsection{Main theorems}

The following simple but important result was one of the first 
discoveries in representation theory. It holds for representations of
any group, not necessarily compact. 

\begin{theorem}[Schur Lemma] \label{t:schur}
  Let $(\pi_i,\,V_i),\,i=1,2,$ be
  any two irreducible finite-dimensional representations of the same
  group $G$. Then any intertwiner $A\colon V_1 \to V_2$ is either
  invertible or zero.
\end{theorem}
\begin{corollary}
If $V$ is an irreducible f.d. representation, then any intertwiner
$A\colon V\to V$ is scalar: $A=c\cdot \id$, $c\in \C$. 
\end{corollary}
\begin{corollary}
  Every irreducible representation of a commutative group is
  one-dimensional. 
\end{corollary}
The following theorem is the fundamental results of the
representation theory of compact groups. Its proof is based on the
technique of invariant integrals on a compact group, which will be
discussed in the next section.

\begin{theorem}\label{t:complete_reducibility}
\par\indent
\begin{enumerate}
\item 
 Any f.d. representation of a compact group is 
 equivalent to a unitary representation. 
\item  Any f.d. representation is completely reducible: it
  can be decomposed into direct sum
  $$
  V=\bigoplus n_i V_i
  $$
  where $V_i$ are pairwise non-equivalent unirreps. Numbers
  $n_i\in \Z_+$  are called {\bf multiplicities}.
\end{enumerate}
\end{theorem}

\subsection{Examples of representations}\label{s:examples_reps}

The representation theory looks rather different for abelian
(i.e. commutative) and non-abelian groups. Here we consider two
simplest examples of both kinds.

Our first example is a 1-dimensional compact connected Lie group. 
Topologically it is a circle which we realize as a set $\T\simeq
\U(1)$ of all complex numbers $t$ with absolute value 1.

Every unirrep of $\T$ is one-dimensional; thus, it is just a
continuous multiplicative map $\pi$ of
$\T$ to itself. It is well-known that every such a map has the form
$$
\pi_k(t)=t^k\quad \text{for some} \quad k\in \Z.
$$ 
The collection of all unirreps of $\T$ is itself a group, called
{\bf Pontrjagin dual} of $\T$ and denoted by $\widehat{\T}$. This
group is isomorphic to $\Z$.
  
By \thref{t:complete_reducibility}, any f.d.
representation $\pi$ of $\T$ is equivalent to a direct sum of
1-dimensional unirreps. So, an equivalence class of $\pi$ is defined
by the multiplicity function $\mu$ on $\widehat{\T} =\Z$ taking
non-negative values:
$$
\pi\simeq \sum_{k\in \Z} \mu(k)\cdot \pi_k.
$$

The many-dimensional case of compact connected abelian Lie group can
be treated in a similar way. Let $T$ be a torus, i.e. an abelian
compact group, $\liet=\Lie(T)$. Then every irreducible representation
of $T$ is one-dimensional and thus is defined by a group homomorphism
$\chi\colon T\to \T^1=\U(1)$. Such homomorphisms are called {\bf
characters}  of $T$. One easily sees that such characters themselves
form a group (Pontrjagin dual of $T$). If we denote by $L$ the kernel
of the exponential map $\liet\to T$ (see \leref{l:L}), one easily sees
that every character has a form 
$$
\chi(\exp(t))=e^{\<t,\la\>}\qquad \la\in i\liet^*, 
\quad \la\in X(T) 
$$
where $X(T)\subset \liet^*$ is the lattice defined by \eqref{e:X}.
Thus,  we can identify the group of characters $\widehat T$  with
$X(T)$. In particular, this shows that $\widehat T\simeq \Z^{\dim
T}$. 

\medskip

The second example is the group $G=\SU(2)$, the simplest
connected, simply-connected non-abelian compact Lie group.
Topologically $G$ is a 3-dimensional sphere since the general element
of $G$ is a matrix of the form
$$
g=\begin{pmatrix}       a &b\\ 
             -\overline b &\overline a
  \end{pmatrix},\quad 
  a,\,b\in \C, \quad |a|^2+|b|^2=1.
$$ 
Let $V$ be 2-dimensional complex vector space, realized by column
vectors $\begin{pmatrix} u\\v\end{pmatrix}$. The group $G$ acts
naturally on $V$. 
This action induces the representation $\Pi$ of $G$ in the space
$S(V)$ of all polynomials in $u,\,v$. It is infinite-dimensional, but 
has many f.d. subrepresentations. In particular, let $S^k(V)$, or
simply $S^k$, be the space of all homogeneous polynomials of degree
$k$. Clearly,\ $\dim S^k = k+1$.

It turns out that the corresponding f.d. representations
$(\Pi_k,\,S^k)$, $k\ge 0$, are irreducible, pairwise non-equivalent
and exhaust the set $\widehat G$ of all unirreps. 

Some particular case are of special interest:

\begin{enumerate}
\item  $k=0$. The space $V_0$ consists of constant functions and
  $\Pi_0$ is the trivial 1-dimensional representation: $\Pi_0(g)\equiv
  1$.

\item  $k=1$. The space $V_1$ is identical to $V$ and $\Pi_1$ is just
  the tautological representation $\pi(g)\equiv g$.


 \item $k=2$. The space $V_2$ is spanned by monomials
  $u^2,\,uv,\,v^2$. 
The remarkable fact is that this representation is equivalent to a
real one. Namely, in the new basis $x=\frac {u^2+v^2}2,\,y=\frac
{u^2-v^2}{2i},\,z= iuv$ we have 
$$
\Pi_2\begin{pmatrix} a&b\\ -\overline b&\overline
a\end{pmatrix}=\begin{pmatrix}
\Re (a^2+b^2)&2\Im (ab)& \Im(b^2-a^2)\\ 2\Im(a\overline b)&
|a|^2-|b|^2& 2\Re (a\overline b)\\ \Im(a^2+b^2)& 2\Re
(ab)&\Re(a^2-b^2)
\end{pmatrix}.
$$
This formula defines a homomorphism $\Pi_2\colon \SU(2)\to \SO(3)$.
It can be shown that this homomorphism is surjective, and its kernel
is the subgroup $\{\pm 1\}\subset \SU(2)$:
$$
1\to \{\pm 1\}\hookrightarrow \SU(2)\overset{\Pi_2}{\longrightarrow} 
\SO(3)\to 1.
$$
The simplest way to see it is to establish the equivalence of $\Pi_2$ 
with the adjoint representation of $G$ in $\g$. The corresponding 
intertwiner is  
$$
S^2\ni(\alpha+i\gamma) u^2+2i\beta
uv+(\alpha-i\gamma)v^2\longleftrightarrow 
\begin{pmatrix} i\beta& \alpha+i\gamma\\ -\alpha+i\gamma&
-i\beta\end{pmatrix}
\in \g.
$$ 
Note that $\SU(2)$ and $\SO(3)$ are the only compact groups associated 
with the Lie algebra $\sll(2,\,\C)$. 

\end{enumerate}

The group $G$ contains the subgroup $H$ of diagonal matrices,
isomorphic
to $\T^1$. Consider the restriction of $\Pi_n$ to $\T^1$. It splits
into the sum of unirreps $\pi_k$ as follows:
$$
\text{Res}\,^G_{\T^1}\Pi_n=\sum_{s=0}^{s=[n/2]} \pi_{n-2s}.
$$
The characters $\pi_k$ which enter  this decomposition are called 
the weights of $\Pi_n$. The collection of all weights (together with 
multiplicities) forms a multiset in $\widehat{\T}$ denoted by
$P(\Pi_n)$, or $P(S^n)$.  

Note the following features of this multiset:

\begin{enumerate}
  \item $P(\Pi_n)$ is invariant under reflection $k\mapsto -k$.

  \item All weights of $\Pi_n$ are congruent modulo 2.

  \item The non-equivalent unirreps have different multisets of
    weights.
\end{enumerate}

Below we show how these features are generalized to all compact
connected Lie groups.

\section{Fourier transform} \label{s:fourier}
\subsection{Haar measure and invariant integral}
The important feature of compact groups is the existence of
so-called ``invariant integral'', or ``average''.

\begin{theorem}\label{t:haar}
  For every compact Lie group $G$, there exists a
  unique measure $dg$  on $G$, called {\bf Haar measure},  which is
  invariant under left  shifts $L_g\colon h\mapsto gh$ and satisfies
  $\int_G dg=1$. 

  In  addition, this measure is also invariant under right shifts
  $h\mapsto hg$ and under involution $h\mapsto h^{-1}$.
\end{theorem}

 Invariance of the Haar measure implies that for every integrable
  function $f(g)$, we have
  $$
  \int_G f(g)\, dg = \int_G f(hg)\, dg = \int_G f(gh)\, dg
   =\int_G f(g^{-1})\, dg
  $$

For a finite group $G$ the integral with respect to the Haar measure 
is just averaging over the group:
$$
\int_G f(g)\, dg=\frac{1}{|G|}\sum_{g\in G}f(g)
$$
 For compact connected Lie groups the Haar measure is given by a
 differential form of top degree which is invariant under right and
 left translations.

For a torus $T^n=\R^n/\Z^n$ with real coordinates $\theta_k\in \R/\Z$
or complex coordinates  $t_k= e^{2\pi i\theta_k}$, the Haar measure
is 
$d^n\theta:= d\theta_1d\theta_2\cdots d\theta_n$ or $d^nt:=
\prod_{k=1}^n\frac {dt_k}{2\pi it_k}$.

In particular, consider a central function $f$ (see \thref{t:T/W}). 
Since every conjugacy class contains elements of the maximal torus $T$
(see \thref{t:tori}), such a function is determined by its values on
$T$, and integral of a central function can be reduced to integration
over $T$. The resulting formula is called {\bf Weyl integration
formula}. For $G=\U(n)$ it looks as follows: 
$$
\int_{U(n)} f(g) dg=\frac{1}{n!}\int_T f(t)\prod_{i<j}|t_i-t_j|^2 d^nt
$$
where $T$ is the maximal torus consisting of diagonal matrices
$$
t=\diag(t_1, \dots, t_n), \quad
    t_k=e^{2\pi i\theta_k},
$$
 and $d^nt$ is defined above.

Weyl integration formula for arbitrary compact group $G$ can be
found in \cite{simon} or \cite[Section 18]{bump}.

The main applications of the Haar measure is the proof of
complete reducibility theorem (\thref{t:complete_reducibility})  and
orthogonality relations (see below).
\subsection{Orthogonality relations and Peter-Weyl theorem}

Let $V_1, V_2$ be unirreps of a compact group $G$. 
Taking any linear operator $A\colon V_1\to V_2$
and averaging the expression $A(g):=\pi_2(g^{-1})\circ A\circ\pi_1(g)$
over $G$, we get an intertwining operator $\langle A\rangle=\int_G 
A(g)dg$. Comparing this fact with the Schur lemma, one obtains the 
following fundamental results.

Let $(\pi,\,V)$ be any unirrep of a compact group $G$. Choose any 
orthonormal basis $\{v_k,\,1\leq k\leq \dim V\}$ in $V$ and denote by 
$t_{kl}^V$, or $t^\pi_{kl}$,  the function on $G$ defined by
$$
t_{kl}^V(g)=(\pi(g)v_l,\,v_k).
$$
The functions $t_{kl}^V$ are called {\bf matrix elements} of the
unirrep $(\pi,\,V)$. 
 
\begin{theorem}[Orthogonality relations]\label{t:orthogonality}
\par\indent
\begin{enumerate}
  \item The matrix elements $t_{kl}^V$ are pairwise
    orthogonal and have norm $(\dim\,V)^{-1/2}$ in $L^2(G,\,dg)$.
  \item The matrix elements corresponding to equivalent unirreps
    span the same subspace in $L^2(G,\,dg)$.
  \item  The matrix elements of two non-equivalent unirreps are
    orthogonal.
    
  \item The linear span of all matrix elements of all unirreps is
  dense  in $C(G), C^\infty(G)$ and in $L^2(G,\,dg)$.
  \textup{(}Generalized  Peter-Weyl  theorem\textup{)}.
\end{enumerate}
\end{theorem}
In particular, this theorem implies that the set $\widehat{G}$ of 
equivalence classes of unirreps is countable.

For a f.-d. representation $(\pi,\,V)$ we introduce the {\bf character} 
of $\pi$ as a function 
\begin{equation}\label{e:character}
\chi_\pi(g)=\tr \pi(g)= \sum_{k=1}^{\dim V} t_{kk}^{\pi}(g).
\end{equation}
It is obviously a central function on $G$. 

\begin{remark}Traditionally, in representation theory the word 
``character'' has two different meanings: 1. A  multiplicative map 
from a group to $\U(1)$. 2. The trace of a representation operator 
$\pi(g)$. For 1-dimensional representations both notions coincide.
\end{remark}

>From orthogonality relations we get the following result. 
\begin{corollary*} 
  The characters of unirreps of $G$ form an orthonormal basis in 
the subspace of central functions in $L^2(G,\,dg)$. 
\end{corollary*}

\subsection{Non-commutative Fourier transform}

The non-commutative Fourier transform on a compact group $G$ is
defined as follows. Let $\widehat G$ denote the set of equivalence
classes of unirreps of $G$. Choose for any $\lambda \in \widehat G$ a
representation $(\pi_\lambda,\,V_\lambda)$ of class $\lambda$ and an
orthonormal basis in $V_\lambda$. Denote by $d(\lambda)$ the dimension
of $V_\lambda$. 

We introduce the Hilbert space $L^2(\widehat G)$ as the space of 
matrix-valued functions on $\widehat G$ whose value at a point 
$\lambda\in \widehat G$ belongs to $\Mat_{d(\lambda)}(\C)$. The norm 
is defined as
$$
||F||^2_{L^2(\widehat G)}=\sum_{\lambda\in\widehat G}d(\lambda)\cdot
\tr\,\bigl(F(\lambda)F(\lambda)^*\bigr).
$$

For a function $f$ on $G$ define its Fourier transform $\widetilde f$
as a matrix valued function on $\widehat G$:
$$
\widetilde f(\lambda)=\int_G f(g^{-1})\pi_\lambda(g)dg.
$$ 

Note that in the case $G=\T^1$ this transform associates to a function
$f$ the set of its Fourier coefficients. 
 
In general this transform keeps some important features of Fourier 
coefficients. 

\begin{theorem}
\par\indent
\begin{enumerate}
  \item  For a function $f\in L^1(G,\,dg)$ the Fourier
  transform $\widetilde f$ is well defined and bounded \textup{(}by
  matrix  norm\textup{)} function on $\widehat G$.

  \item For a function $f\in L^1(G,\,dg)\cap L^2(G,\,dg)$ the
    following  analog of the Plancherel formula holds:
  $$
  ||f||^2_{L^2(G,\,dg)}:=\int_G|f(g)|^2dg= \sum_{\lambda\in\widehat G}
  d(\lambda)\cdot \tr\,\big (\widetilde f(\lambda)\widetilde 
  f(\lambda)^*\big )=:||\widetilde f||^2_{L^2(\widehat G)}.
  $$

  \item  The following inversion formula expresses $f$ in terms of
  $\widetilde f$:
  $$
  f(g)=\sum_{\lambda\in\widehat G}d(\lambda)\cdot 
    \tr\,\bigl (\widetilde f(\lambda)\pi_\lambda(g)\bigr ).
  $$

  \item  The Fourier transform send the convolution to the matrix
  multiplication:
  $$
  \widetilde{f_1 * f_2}=\widetilde f_1 \cdot \widetilde f_2
  $$
  where the convolution product $*$ is defined by 
  $$
  (f_1 * f_2)(h)=\int_G f_1(hg)f_2(g^{-1})\, dg. 
  $$
  
\end{enumerate}
\end{theorem}
Note the special case of the inversion formula for $g=e$:
$$
f(e)=\sum_{\lambda\in\widehat G}d(\lambda)\cdot \tr\,\big (\widetilde
f(\lambda)\big ),\quad \text{or}\quad
\delta(g)=\sum_{\lambda\in\widehat G}
d(\lambda)\cdot \chi_\lambda(g)
$$
where $\de(g)$ is Dirac's delta function: $\int_G f(g)\de(g)\,
dg=f(e)$. Thus we get a presentation of Dirac's delta function as a
linear combination of characters.

\section{Classification of finite-dimensional representations}
\label{s:unirreps}
In this section, we give a classification of unirreps of a connected
compact Lie group $G$. 

\subsection{Weight decomposition}
Let $G$ be a connected compact group with maximal torus $T$, and
let $(\pi,\,V)$ be a f.d. representation of $G$. Restricting it to $T$ 
and using complete reducibility, we get the following result.
\begin{theorem}
The vector space $V$ can be written in the form
\begin{equation}\label{e:weight_decompos}
V=\bigoplus_{\la\in X(T)} V_\la,\qquad
V_\la=\{v\in V\st \pi_*(t)v =\<\la,t\>v\quad \forall t \in \liet\}.
\end{equation}
where $X(T)$ is the character group of $T$ defined by \eqref{e:X}. 

The spaces $V_\la$ are called {\bf weight subspaces}, vectors $v\in
V_\la$ --- {\bf weight vectors} of weight $\la$. The set
\begin{equation}
P(V)=\{\la\in X(T)\st V_\la\ne \{0\}\}
\end{equation}
is called the {\bf set of weights} of $\pi$, or the {\bf spectrum} of
$\Res^G_T \pi$, and 
$$
\mult_{(\pi,\,V)}(\la):=\dim V_\la
$$
is called the {\bf multiplicity} of $\la$ in $V$. 
\end{theorem}

The next theorem easily follows from the definition of the Weyl
group. 
\begin{theorem}
For any f. d. representation $V$ of $G$, the set of
weights with multiplicities is invariant under the action of the
Weyl group:
$$
w(P(V))=P(V),\qquad \mult_{(\pi,\,V)}(\la)=\mult_{(\pi,\,V)}(w(\la))
$$
for any $w\in W$. 
\end{theorem}

\subsection{Classification of unirreps}

Recall that $R$ is the root system of $\g_\C$. Assume
that we have chosen a basis of simple roots $\al_1,\dots, \al_r\subset
R$. Then $R=R_+\cup R_-$; roots $\al\in R_+$ can be written as a linear
combination of simple roots with positive coefficients, and
$R_-=-R_+$. 

A (not necessarily finite-dimensional) representation of $\g_\C$ is
called a {\bf highest weight representation} if it is generated by a
single vector $v\in V_\la$ (the highest weight vector)  such that
$\g_\al v=0$ for all positive roots $\al\in R_+$. 

It can be shown  that for every $\la\in X(T)$, there is a unique
irreducible highest weight representation of $\g_\C$ with highest
weight $\la$, which is denoted $L(\la)$. However, this representation
can be infinite-dimensional; moreover, it may not be possible to lift
it to a representation of $G$. 

\begin{definition}
A weight $\la\in X(T)$ is called {\bf dominant} if $\<\la,
\al^\vee_i\>\in\Z_+$ for any simple root $\al_i$. The set of all
dominant weights is denoted by $X_+(T)$.  
\end{definition}

\begin{theorem}
  \par\indent
  \begin{enumerate}
  \item All weights of $L_\la$ are of the form $\mu=\la-\sum
    n_i\al_i$, $n_i\in \Z_+$. 
  \item  Let $\la\in X_+$. Then the irreducible highest weight representation
    $L(\la)$ is finite-dimensional and lifts to a representation of $G$. 
  \item Every irreducible finite-dimensional representation of $G$ is of the
    form $L(\la)$ for some $\la\in X_+$. 
    \end{enumerate}
\end{theorem}
Thus, we have a bijection $\{\text{unirreps of }G\}\leftrightarrow
X_+$. 
\begin{example}
Let $G=\SU(2)$. There is a unique simple root $\al$ and the unique
fundamental weight $\om$, related by $\al=2\om$. Therefore,  $X_+=\Z_+
\cdot \om$ and unirreps are indexed by non-negative integers. The 
representation with highest weight $k\cdot \om$ is precisely the 
representation $\Pi_k$ constructed in \seref{s:examples_reps}. 
\end{example}

\begin{example}
Let $G=\U(n)$. Then $X=\Z^n$, and $X_+=\{(\la_1,\dots, \la_n)\in
\Z^n\st
 \la_1\ge\dots\ge \la_n\}$. Such objects
are well known in combinatorics: if we additionally assume that
$\la_n\ge 0$, then such dominant weights are in bijection with
partitions with $n$ parts. They  can  also be described by {\bf Young
diagrams} with $n$  rows (see \cite{fulton-harris}). 
\end{example}

\subsection{Explicit construction of representations}
In addition to description of unirreps as highest weight representations, they
can also be constructed in other ways. In particular, they can be
defined analytically as follows. Let $B=HN_+$ be the Borel subgroup in
$G_\C$; here $H=\exp \h, N_+=\exp\sum_{\al \in R_+} c_\al e_\al$.
For $\la\in \h^*$, let
$\chi_\la\colon B\to \C^\times$ be a multiplicative map defined by
\begin{equation}
\chi_\la(hn)=e^{\<h,\la\>}.
\end{equation}

\begin{theorem}[E.Cartan--Borel--Weil]
Let $\la\in X(T)$. Denote by $V(\la)$ be the space of complex-analytic
functions on $G_\C$ which satisfy the following transformation
property: 
$$
f(gb)=\chi^{-1}_\la(b)f(g),\quad  g\in G_\C,\quad b\in B.
$$
The group $G_\C$ acts on $V(\la)$ by left shifts:
\begin{equation}
\bigl(\pi(g)f\bigr)(h)=f(g^{-1}h).
\end{equation}

  \begin{enumerate}
  \item $V(\la)\neq \{0\}$ iff $-\la\in X_+$.

  \item If $-\la\in X_+$, the representation of $G$ in $V(\la)$ is
  equivalent to $L(w_0(\la))$, where $w_0\in W$ is the unique element
  of the Weyl group which sends $R_+$ to $R_-$. 
\end{enumerate}
\end{theorem}

This theorem can also be reformulated in more geometric terms: the
spaces $V(\la)$ are naturally interpreted as spaces of global sections
of appropriate line bundles on the {\bf flag variety}
$\mathcal{B}=G_\C/B=G/T$. 

For classical groups irreducible representations can also be
constructed explicitly as the subspaces in tensor powers $(\C^n)^{\ttt
k}$,  transforming in a certain way under the action of the symmetric
group $S_k$.

\section{Characters and multiplicities}
\label{s:multiplicities}
\subsection{Characters}

Let $(\pi,V)$ be a finite-dimensional representation of $G$ and let $\chi_\pi$
be its  character as defined by \eqref{e:character}. Since $\chi_\pi$ is
central, and every element in $G$ is conjugate to an element of $T$, 
$\chi_\pi$ is completely determined by its restriction to $T$, which can be
computed from the weight decomposition \eqref{e:weight_decompos}: 
\begin{equation}\label{e:character2}
\chi_\pi|_T=\sum_{\la\in X(T)}\dim V_\la\cdot e_\la = \sum_{\la\in 
X(T)}\mult_\pi \la\cdot e_\la
\end{equation}
where $e_\la$ is the function on $T$ defined by $e_\la(\exp(t))=e^{\<t,\la\>}$,
$t\in \liet$. Note that  $e_{\la+\mu}=e_\la e_\mu$ and that
$e_0=1$. 

\subsection{Weyl character formula}
\begin{theorem}[Weyl character formula]
Let $\la\in X_+$. Then 
$$
\chi_{L_\la}=\frac{A_{\la+\rho}}{A_\rho}, \qquad 
A_\mu=\sum_{w\in W} \varepsilon(w)e_{w(\mu)}
$$
where, for $w\in W$, we denote  $\varepsilon(w)=\det w$ considered as
a linear map $\liet^*\to\liet^*$, and $\rho=\frac{1}{2}\sum_{R_+}\al$.
\end{theorem}
In particular, computing the value of the character at point $t=0$ by
L'Hopital's rule, it is possible to deduce  the following formula
for the dimension of irreducible representations:
\begin{equation}
\dim L(\la)=\prod_{\al\in
R_+}\frac{\<\al^\vee,\la+\rho\>}{\<\al^\vee,\rho\>}. 
\end{equation}

\begin{example}
Let $G=\SU(2)$. Then Weyl character formula gives, for irreducible
representation $\Pi_k$ with highest weight $k\cdot\om$
$$
\chi_{\Pi_k}=\frac{x^{k+1}-x^{-(k+1)}}{x-x^{-1}}
=x^k+x^{k-2}+\dots +x^{-k}, \qquad 
x=e_{\om}
$$
which implies $\dim \Pi_k=k+1$. 
\end{example}

Weyl character formula is equivalent to the following formula for
weight multiplicities, due to Kostant:
$$
\mult_{L(\la)}\mu=\sum_{w\in W}\epsilon(w)K(w(\la+\rho)-\rho-\mu),
$$
where $K$ is Kostant's partition function: $K(\tau)$ is the number of
ways of writing $\tau$ as a sum of positive roots (with repetitions).

For classical Lie groups such as $G=\U(n)$, there are more explicit
combinatorial formulas for weight multiplicities; for $\U(n)$, the
answer can be written in terms of the number of {\bf Young tableaux}
of a given shape. Details can be found in \cite{fulton-harris}.

\subsection{Tensor product multiplicities}
Let $(\pi,V)$ be a finite-dimensional representation of $G$. By
complete reducibility, one can write $V=\sum n_\la L(\la)$. The
coefficients $n_\la$ are called multiplicities; finding them is an
important problem in many applications. In particular, a special case
of this is finding the multiplicities in tensor product of two
unirreps:
$$
L(\la)\ttt L(\mu)=\sum N_{\la\mu}^\nu L(\nu). 
$$

Characters provide a practical tool of computing multiplicities:
since characters of unirreps are linearly independent, multiplicities
can be found from the condition that $\chi_V=\sum n_\la
\chi_{L(\la)}$. In particular, 
$$
\chi_{L(\la)}\chi_{L(\mu)}=\sum N_{\la\mu}^\nu \chi_{L(\nu)}.
$$
\begin{example}
For $G=\SU(2)$, tensor product multiplicities are given by 
$$
\Pi_n\ttt \Pi_m=\bigoplus \Pi_l
$$
where the sum is taken over all $l$ such that $|m-n|\le l\le m+n$,
$m+n+l$ is even. 
\end{example}
For $G=\U(n)$, there is an algorithm for finding the tensor product
multiplicities, formulated in the language of Young tableaux
(Littlewood--Richardson rule).  There are also tables and computer
programs for computing these multiplicities; some of them are listed
in the bibliography. 


\section*{See Also}
Classical groups and homogeneous spaces. Lie groups: general theory. 
Infinite-dimensional Lie algebras.

\section*{Keywords}
Lie groups; compact groups; maximal tori; Lie algebras; characters;
representations; roots; weights;  highest weight; multiplicities; 

\bibliographystyle{amsalpha}

\end{document}